\documentclass[9pt,amstex]{article}
\usepackage{amssymb,amsmath} 
\usepackage{latexsym}

\newcommand{\ds}{\displaystyle}
\newcommand{\cbi}[2]{\left(\begin{array}{c}#2\\#1\end{array}\right)}

\newtheorem{dfn}{Definition}[section] 
\newtheorem{rmk}{Remark}[section]
\newtheorem{thm}{Theorem}[section] 
\newtheorem{cor}{Corollary}[section]
\newtheorem{prop}{Proposition}[section] 
\newtheorem{lem}{Lemma}[section]

\def \als{\mathop{\rightharpoonup}\limits}

\def \ars{\mathop{\leftharpoonup}\limits}
\def\buildrel#1_#2^#3{\mathrel{\mathop{\kern 0pt#1}\limits_{#2}^{#3}}}

\newcommand{\Pf}{{\em Proof}. }
\newcommand{\EPf}
{%
\mbox{}%
\nolinebreak%
\hfill%
\rule{2mm}{2mm}%
\medbreak%
\par%
}

\newcommand{\auto}{\mbox{$\mathtt{Aut}$}}
\newcommand{\End}{\mbox{$\mathtt{End}$}}

\newcommand{\der}{\mbox{$\mathtt{Der}$}}
\newcommand{\id}{\mbox{$\mathtt{Id}$}}
\newcommand{\Aff}{\mbox{$\mathtt{Aff}$}}
\newcommand{\Symp}{\mbox{$\mathtt{Symp}$}}

\newcommand{\C}{\mathbb C} 

\newcommand{\R}{\mathbb R}

\newcommand{\g}{{\mathfrak{g}}{}} 
 
\renewcommand{\l}{{\mathfrak{l}}{}}
\renewcommand{\k}{{\mathfrak{k}}{}} 
\newcommand{\hg}{{\mathfrak{h}}({\mathfrak{g}}){}}
\newcommand{\p}{{\mathfrak{p}}{}} 
\newcommand{\q}{{\mathfrak{q}}{}}

\newcommand{\s}{{\mathfrak{s}}{}}

\renewcommand{\d}{{\mathfrak{d}}{}} 
\newcommand{\h}{{\mathfrak{h}}{}} 
\renewcommand{\a}{{\mathfrak{a}}{}} 
\renewcommand{\b}{{\mathfrak{b}}{}} 

\newcommand{\z}{{\mathfrak{z}}{}}

\newcommand{\pw}{{\bf w}{}}
\newcommand{\EK}{{\mathfrak{l}}{}}
\newcommand{\KK}{{\mathfrak{k}}{}} 

\def\tref#1{Theorem~\ref{#1}}

\def\cref#1{Corollary~\ref{#1}}
\def\dref#1{Definition~\ref{#1}}

\title{Universal Deformation Formulae, Symplectic Lie groups and
Symmetric Spaces}

\author{{\bf Pierre Bieliavsky}\\
Universit\'e Libre de Bruxelles, Belgium,\\
pbiel@ulb.ac.be\\
{\bf Philippe Bonneau}\\
Universit\'e de Bourgogne, France,\\
bonneau@u-bourgogne.fr\\
{\bf Yoshiaki Maeda}\\
Keio University, Japan,\\
maeda@math.keio.ac.jp}
\date{June 25, 2003} 

\begin{document}
\maketitle

\begin{abstract}
We apply methods from strict quantization of solvable
symmetric spaces to obtain universal deformation formulae for actions
of a class of solvable Lie groups.  We also study compatible 
co-products by generalizing the notion of smash product in the context
of Hopf algebras.
\end{abstract}

\section{Introduction}
Let $G$ be a group acting on a set $M$. Denote by $\tau:G\times M\to 
M:(g,x)\mapsto\tau_g(x)$ the (left) action and by
$\alpha:G\times\mbox{Fun}(M)\to\mbox{Fun}(M)$ the corresponding action
on the space of (complex valued) functions (or formal series) on $M$
($\alpha_g:=\tau^\star_{g^{-1}}$).  Assume that on a subspace
${\mathbb A}\subset\mbox{Fun}(G)$, one has an associative $\C$-algebra
product $\star^G_{\mathbb A}:{\mathbb A}\times {\mathbb
A}\to {\mathbb A}$ such that
\begin{enumerate}
\item[(i)] ${\mathbb A}$ is invariant under the (left) regular action of $G$ on 
$\mbox{Fun}(G)$,
\item[(ii)] the product $\star^G_{\mathbb A}$ is left-invariant as well
i.e. for all $g\in G; a,b\in {\mathbb A}$, one has
\begin{equation}
(L_g^\star a) \star^G_{\mathbb A} (L_g^\star b)=L_g^\star(a\star^G_{\mathbb A}b).
\end{equation}
\end{enumerate}

\noindent Given a function on $M$, $u\in\mbox{Fun}(M)$, and a point $x\in M$, one 
denotes by $\alpha^x(u)\in\mbox{Fun}(G)$ the function on $G$ defined as 
\begin{equation}
\alpha^x(u)(g):=\alpha_g(u)(x).
\end{equation}
Then one readily observes that the subspace ${\mathbb B}\subset\mbox{Fun}(M)$ defined as 
\begin{equation}
{\mathbb B}:=\{u\in\mbox{Fun}(M)\, |\, \forall x\in M:\alpha^x(u)\in
{\mathbb A}\}
\end{equation}
becomes an associative $\C$-algebra when endowed
with the product $\star^M_{\mathbb B}$ given by
\begin{equation}\label{UDF}
u\star^M_{\mathbb B}v(x):=(\alpha^x(u)\star^G_{\mathbb A}\alpha^x(v))(e)
\end{equation}
($e$ denotes the neutral element of $G$).
Of course, all this can be defined for right actions as well.
\begin{dfn}
Such a pair $({\mathbb A},\star^G_{\mathbb A})$ is called a (left) {\bf universal deformation}
of $G$, while Formula (\ref{UDF}) is called the associated {\bf universal 
deformation formula} (briefly {\bf UDF}).
\end{dfn}
In the present article, we will be concerned with the case where $G$ is a Lie 
group. The function space ${\mathbb A}$ will be either 

\noindent - a functional subspace (or a topological completion) of 
$C^\infty(G,\C)$ containing the smooth compactly supported functions in which case we will talk
about {\bf strict deformation} (following Rieffel \cite{Rieffel1}), 

\noindent or, 

\noindent - the space ${\mathbb A}=C^\infty(G)[[\hbar]]$ of formal power series with coefficients 
in the smooth functions on $G$ in which case, we'll speak about {\bf 
formal deformation}. In any case, we'll assume the product $\star^G_{\mathbb A}$ 
admits an asymptotic expansion of  star-product type:
\begin{equation*}
a\star^G_{\mathbb A}b\sim
ab+\frac{\hbar}{2i}\pw({\rm d}u,{\rm d}v)+o(\hbar^2)\qquad (a,b\in C^\infty_c(G)),
\end{equation*}
where $\pw$ denotes some (left-invariant) Poisson bivector on
$G$\cite{BFFLS}.  In the strict cases considered here, the product
will be defined by an integral three-point kernel $K\in
C^\infty(G\times G\times G)$:
\begin{equation*}
a\star^G_{\mathbb A}b(g):=\int_{G\times G}a(g_1)\,b(g_2)K(g,g_1,g_2){\rm d}g_1\,{\rm 
d}g_2\qquad(a,b\in {\mathbb A})
\end{equation*}
where ${\rm d}g$ denotes a normalized left-invariant Haar 
measure on $G$. Moreover, our kernels will be
of  {\bf WKB type} \cite{Wein, Kara} i.e.:
\begin{equation*}
K=A\,e^{\frac{i}{\hbar}\Phi},
\end{equation*}
with $A$ (the {\bf amplitude}) and $\Phi$ (the {\bf phase}) in $C^\infty(G\times 
G\times G,\R)$ being invariant under the (diagonal) action by 
left-translations.

\noindent Note that in the case where the group $G$  acts smoothly on a smooth
manifold $M$ by diffeomorphisms: $\tau:G\times M\to 
M:(g,x)\mapsto\tau_g(x)$, the first-order expansion term of
$u\star^M_{\mathbb B}v,\quad u,v\in C^{\infty}(M)$ defines a Poisson
structure $\pw^{M}$ on $M$ which can be expressed in terms of a basis $\{X_i\}$ of 
the Lie  algebra $\g$ of $G$ as:
\begin{equation}\label{PS}
    \pw^{M}=\left[\pw_e\right]^{ij}X^\star_i\wedge
X^\star_j,
\end{equation}
where $X^\star$ denotes the fundamental vector field on $M$
associated to $X\in\g$.

\vspace{5mm}

\noindent Strict deformation theory in the WKB context was initiated 
by Rieffel in \cite{Rieffel2} in the cases where $G$ is either Abelian or
1-step nilpotent.  Rieffel's work has led to what is now called 
`Rieffel's machinery'; producing a whole class of exciting non-commutative
manifolds (in Connes sense) from the data of Abelian group actions on
$C^{\star}$-algebras \cite{CoLa}.

\noindent The study of formal UDF's for non-Abelian group actions in
our context was initiated in \cite{GiaZha}
where the case of the group of  affine transformations of the real line
(`$ax+b$') was explicitly described.

\noindent In the strict (non-formal) setting, UDF's for Iwasawa subgroups
of $SU(1,n)$ have been explicitly given in \cite{BiMas}.  These were
obtained by adapting a method developed by one of us in the symmetric space
framework \cite{Biep00a}.

\noindent This work has three parts.  We first give some results about the
structure of symplectic symmetric spaces for which the action of the
linear holonomy algebra at a point  $o$ has  isotropic range in
the symplectic vector space tangent at $o$.  In particular, we show
that, when split (cf.  Definition \ref{SPLIT}), such a space can be
identified with the underlying manifold of a Lie subgroup $S$ of its
automorphism group.  In particular, this yields a distinguished class
of symplectic Lie groups in Lichn\'erowicz sense which we call {\sl
elementary}.

\noindent Secondly, we use the strict quantization method of
\cite{Biep00a} to deduce universal deformation formulae for such
symplectic Lie groups.  We extend these results to a class of Abelian
extensions of elementary symplectic Lie groups.  We therefore obtain
UDF's for a class of groups that can be thought of generalized $ax+b$
groups.  

\noindent Finally, we generalize the classical definition of smash products
in the context of (formal) bialgebras (or Hopf algebras) by defining a
class of products constructed from the data of {\sl bi}-module algebras;
we call these products {\sl L-R smash} products.  We formally realize
each of our UDF's as a L-R-smash product.  As a corollary, we obtain 
compatible co-products.

\section{General facts about symplectic symmetric spaces}\label{SSS}

\begin{dfn}\label{SSSDEF}\cite{Biep95a} A {\bf symplectic symmetric space} 
is a triple $(M, \omega,s)$, where $(M,\omega)$ is a smooth connected 
symplectic manifold, and where $s : M \times M \to M$ is a smooth map such 
that 
  
\begin{enumerate} 
\item[(i)] for all $x$ in $M$, the partial map $s_x : M \to M : y \mapsto 
s_x (y) := s(x,y)$ is an involutive symplectic diffeomorphism of $(M,\omega)$ 
called the {\bf symmetry} at $x$. 
\item[(ii)] For all $x$ in $M$, $x$ is an isolated fixed point of $s_x$. 
\item[(iii)] For all $x$ and $y$ in $M$, one has $s_xs_ys_x=s_{s_x (y)}$. 
\end{enumerate}

\end{dfn}

\begin{dfn} Two symplectic symmetric spaces  $(M,\omega,s)$ 
and $(M',\omega ',s')$ are {\bf isomorphic} if there exists a symplectic 
diffeomorphism $\varphi: (M,\omega) \rightarrow (M',\omega')$ such that 
$\varphi s_x=s'_{\varphi (x)} \varphi$. Such a $\varphi $ is called an {\bf
isomorphism} of $(M,\omega,s)$ onto $(M',\omega',s')$. When 
$(M,\omega,s) = (M',\omega',s')$, one talks about {\bf automorphisms}. 
The group of all automorphisms of the symplectic symmetric space 
$(M,\omega,s)$ is denoted by $\auto(M,\omega,s)$.
\end{dfn}

\begin{prop} On a symplectic symmetric space $(M,\omega,s)$, there exists 
a unique affine connection $\nabla$ which is invariant under the
symmetries.  Moreover, this connection satisfies the following
properties.

\begin{enumerate} 
\item[(i)] For all smooth tangent vector fields $X,Y,Z$  on $M$ and all 
points $x$ in $M$, one has $$ \omega_{x} (\nabla _X Y,Z)=\frac{1}{2}
X_x.\omega (Y+s_{x_{\star}}Y,Z).$$ \item[(ii)] $(M,\nabla )$ is an affine 
symmetric space. In particular $\nabla$ is torsion free and its curvature 
tensor is parallel.  
\item[(iii)] The symplectic form $\omega$ is parallel; $\nabla$ is therefore 
a symplectic connection.   
\item[(iv)] One has
$$\auto(M,\omega,s)=\auto(M,\omega,\nabla)=\Aff(\nabla)\cap\Symp(\omega).$$
\end {enumerate} 

\end{prop}

The connection $\nabla$ on the symmetric space $(M,s)$ is called the {\bf 
Loos connection}. The following facts are classical (see \cite{Lo}).

\begin{thm}\label{TG}
Let $(M,\omega,s)$ be a symplectic symmetric space and $\nabla$ its 
Loos connection. Fix $o$ in $M$ and denote by $H$ the stabilizer of 
$o$ in $\auto(M,\omega,s)$. Denote by $G$ the {\bf transvection group} of 
$(M,s)$ (i.e.~the subgroup of $\auto(M,\omega,s)$ generated by $\{ s_x \circ 
s_y \, ; \, x,y \in M \}$)  and set $K=G \cap H$. Then, 

\begin{enumerate} 
\item[(i)] the transvection group $G$ is a connected Lie transformation 
group of $M$. It is the smallest subgroup of $\auto(M,\omega,s)$ which is transitive 
on $M$ and stabilized by the conjugation $\tilde{\sigma} : 
\auto(M,\omega,s) \to \auto(M,\omega,s)$ defined by $\tilde{\sigma}(g)=s_o g s_o$.
\item[(ii)] The homogeneous space $M=G/_{\textstyle{K}}$ is reductive. The 
Loos connection $\nabla$ coincides with the canonical connection induced 
by the structure of reductive homogeneous space. 
\item[(iii)] Denoting by $G^{\tilde{\sigma}}$ the set of $\tilde{\sigma}$-fixed 
points in $G$ and by $G_0^{\tilde{\sigma}}$ the connected component of
the identity,
one has $$ G_0^{\tilde{\sigma}} \subset K \subset G^{\tilde{\sigma}}.$$ 
The Lie algebra $\k$ of $K$ is isomorphic to the
holonomy algebra with respect to the canonical connection $\nabla$. 
\item[(iv)] Denote by $\sigma$ the involutive automorphism of the Lie
algebra $\g$ of $G$ induced by the automorphism $\tilde{\sigma}$. 
Denote by $\g = \k \oplus \p$ the decomposition in $\pm 1$-eigenspaces
for $\sigma$.  Then, identifying $\p$ with $T_o(M)$, one has $$ exp
(X) = s_{Exp_o(\frac{1}{2} X)} \circ s_o $$ for all $X$ in a
neighborhood of $0$ in $\p$.  Here $exp$ is the exponential map $exp :
\g \to G$ and $Exp_o$ is the exponential map at point $o$ with respect
to the connection $\nabla$.

\end{enumerate}
\end{thm}

\begin{dfn} Let $(\g,\sigma)$ be an {\bf involutive algebra}, that is, 
$\g$ is a finite dimensional real Lie algebra and $\sigma$ is an
involutive automorphism of $\g$.  Let $\Omega$ be a skewsymmetric
bilinear form on $\g$.  Then the triple $(\g,\sigma, \Omega)$ is
called a {\bf symplectic triple} if the following properties are
satisfied.

\begin{enumerate}	
\item[(i)] Let $\g=\k\oplus\p$ where $\k$ (resp.  $\p$) is the $+1$
(resp.  $-1$) eigenspace of $\sigma$.  Then $[\p,\p]=\k$ and the
representation of $\k$ on $\p$, given by the adjoint
action, is faithful.  \item[(ii)] $\Omega$ is a Chevalley 2-cocycle
for the trivial representation of $\g$ on $\R$ such that for any $X$
in $\k$, $i(X){\Omega}=0$.  Moreover, the restriction of $\Omega$ to
$\p \times\p$ is nondegenerate.
\end{enumerate} 

The dimension of $\p$ defines the {\bf dimension}	of	the triple. 
Two	such triples $({\g}_i,\sigma_i,{\Omega}_i)$ $(i=1,2)$ are 
{\bf isomorphic} if there exists a Lie algebra isomorphism 
$\psi :\g_1\rightarrow\g_2$ such that $\psi \circ \sigma_1 = 
\sigma_2 \circ \psi$ and $\psi^*{\Omega}_2={\Omega}_1$. 
\end{dfn}

\tref{TG} associates to a symplectic symmetric space $(M,\omega,s)$ 
an involutive Lie algebra $(\g,\sigma)$. Denoting by $\pi : G \to M$ 
the natural projection, one checks that the triple $(\g, \sigma, \Omega = 
\pi^* (\omega_o))$ is a symplectic triple. This implies the next 
proposition.  

\begin{prop} There is a bijective correspondence between the isomorphism
classes of simply connected symplectic symmetric spaces $(M,\omega,s)$ 
and the isomorphism classes of symmetric triples $(\g,\sigma,{ 
\Omega})$.  
\end{prop}

Since a symmetric symplectic manifold $(M,\omega,s)$ is a symplectic 
homogeneous space of its transvection group $G$, it seems natural, when 
possible, to relate $(M,\omega,s)$ to a coadjoint orbit of $G$ in ${\cal 
G}^\star$. Recall first the two following definitions.  

\begin{dfn} Let $G$ be a Lie group of symplectomorphisms acting on a 
symplectic manifold $(M,\omega)$. For every element $X$ in the Lie algebra 
$\g$ of $G$, one denotes by $X^*$ the {\bf fundamental vector field} 
associated to $X$, i.e.~for $x$ in $M$,
$$
X^*_x = \frac{d}{dt} \exp (-tX) x|_{t=0}.
$$
The action is called {\bf weakly Hamiltonian} if for all $X$ in $\g$ 
there exists a smooth function $\lambda_X \in C^\infty(M)$ such that 
$$
i(X^*)\omega = d \lambda_X.
$$ 
In this case, if the correspondence $\g \to C^\infty (M) : X \mapsto 
\lambda_X$ is also a homomorphism of Lie algebras, one says that the action 
of $G$ on $(M,\omega)$ is {\bf Hamiltonian}. (The Lie algebra structure 
on $C^\infty (M)$ is defined by the Poisson bracket.)
\end{dfn}

\begin{prop} Let $t=(\g,\sigma,\Omega)$ be a symplectic triple 
and let $(M,\omega, s)$ be the associated simply connected symplectic 
symmetric space. The action of the transvection group $G$ on $M$ is 
Hamiltonian if and only if ${\Omega}$ is a Chevalley coboundary, that 
is, there exists an element $\xi$ in $\g^*$ such that $\Omega = 
\delta \xi$. In this case, $(M,\omega,s)$ is a $G$-equivariant symplectic 
covering of ${\cal O}$, the coadjoint orbit of $\xi$ in $\g^\star$. 
\end{prop}

The action of the transvection group $G$ is in general not Hamiltonian. 
We therefore need to consider a one-dimensional central extension of $G$ 
rather than $G$ itself. At the infinitesimal level, this corresponds to 
extending the algebra $\g$ by the 2-cocycle $\Omega$. This way, one associates 
to any symplectic symmetric space an exact triple in the following sense 
(see \cite{Bi1} for details).

\begin{dfn} An {\bf exact triple} is a triple $\tau = (\h,\sigma,{\bf \Omega})$ where
\begin{itemize} 
\item $(\h,\sigma)$ is an involutive Lie algebra such that if $\h=\EK
\oplus \p$ is the decomposition w.r.t. $\sigma$ one has $\left[\p,\p\right]=\EK,$ 
\item ${\bf \Omega}$ is a Chevalley 2-coboundary 
(i.e. ${\bf \Omega}=\delta \xi \,  , \quad
\xi \in \h^\star$) such that $i(\EK){\bf \Omega}=0$ and ${\bf \Omega}|_{\p \times \p}$
is symplectic. 
\end{itemize} 
\end{dfn} 

\begin{rmk} One can choose $\xi \in \h^\star$
such that $\xi(\p)=0$. 
\end{rmk} 

\begin{lem}\label{ET} Let $t=(\g,\sigma,{\bf\Omega})$ be a symplectic triple. 
Assume that the triple $t$ is non-exact.  Consider the triple
$\tau(t)=(\hg,\sigma_{\hg},{\bf \Omega}_{\hg})$ constructed as follows:
\begin{itemize} 
\item $0\to\R E \to \hg \stackrel{\pi}{\to} \g\to0$ is the central
extension defined by $\left[ X,Y \right]_{\hg} = {\bf \Omega}(X,Y) E
\oplus \left[ X,Y \right]_{\g}$ \item $\sigma_{\hg} = id_{\R E} \oplus
\sigma$ \item ${\bf \Omega}_{\hg}$ is the trivial extension of ${\bf
\Omega}$ to $\hg$.
\end{itemize}
Then, $\tau(t)$ is an exact triple. 
\end{lem} 

\begin{rmk}\label{DIM(Z)}{\rm 
Observe that, when associated to a (transvection) symplectic 
triple, the center $\z(\hg)$ of the Lie algebra $\hg$ occurring in an exact 
triple is at most one dimensional. Indeed, on the one hand, exactness 
implies $\z(\hg)\subset\k$. One the other hand, faithfulness of the holonomy 
representation forces $\dim(\z(\hg)\cap\k)\leq1$ since $\k$ is either the 
holonomy algebra itself or a one dimensional central extension.
}\end{rmk}

\section{Elementary solvable symplectic symmetric spaces and their 
strict quantization}\label{SQ}

In \dref{ESET} below, we define a particular type of solvable 
symmetric spaces which we call elementary.  It has been proven
(\cite{Bi1}, Proposition~3.2) that every solvable symmetric space is realized
through a sequence of split extensions by Abelian (flat) factors
successively taken over an elementary solvable symmetric space.  We
therefore consider elementary solvable symmetric spaces as the ``first
induction step'' when studying solvable symmetric spaces.

\begin{dfn}\label{ESET}
A symplectic symmetric space $(M,\omega,s)$ is called an {\bf elementary 
solvable} symplectic symmetric space if its associated exact triple 
$(\hg,\sigma,\Omega = \delta \xi)$ (see Lemma \ref{ET}) is of the
following type.
\begin{enumerate}
\item[(i)] The Lie algebra $\hg$ is a split extension of Abelian Lie 
algebras $\a$ and $\b$~:
$$
0\to\b\longrightarrow\hg{\longrightarrow}\a\to0.
$$
\item[(ii)] The automorphism $\sigma$ preserves the splitting 
$\hg=\b\oplus\a$.
\end{enumerate}
Such an exact triple (associated to an elementary solvable symplectic 
symmetric space) is called an {\bf elementary solvable exact triple}. 
\end{dfn}

Observe that, since $\a \cap \KK \subset \a \cap [\hg,\hg] = 0$, one 
has $\a \subset \p$. Therefore $\b = \KK \oplus \EK$, with $\EK 
\subset \p$. Moreover, since $\EK$ and $\a$ are Abelian and 
$\Omega$ is nondegenerate, the subspaces $\a$ and $\EK$ of $\p$ are 
dual Lagrangians.

Now let $(M,\omega,s)$ be an elementary solvable symplectic symmetric 
space with associated exact triple $(\hg,\sigma,\Omega = \delta \xi)$
as above. In a neighborhood $U$ of the origin, the map
\begin{equation}\label{COORD}
\p=\a\times\l\to M: (a,l)\mapsto\exp(a)\exp(l).o
\end{equation}
turns out to be a Darboux chart when $U\subset\p$ has the symplectic
structure $\Omega=\delta\xi$.  Moreover, there exists a unique
immersion $\phi:U\cap\a\to\a$ such that in the local coordinate system
(\ref{COORD}), one has the following linearization property:
\begin{equation}\label{TWIST}
\xi(\sinh(a)l)=\xi[\phi(a),l];
\end{equation}
where, for $a\in\a$ we set 
$\sinh(a):=\frac{1}{2}(\exp(\rho(a))-\exp(-\rho(a)))\in\End(\b)$.
This immersion is called the {\bf twisting map}.
\begin{prop}
An elementary solvable symplectic symmetric 
space is strictly geodesically convex if and only if its 
associated twisting map extends to $\a$ as a global diffeomorphism of $\a$.
In this case, the Darboux chart (\ref{COORD}) extends as a global 
symplectomorphism $(\p,\Omega)\to(M,\omega)$.
\end{prop}
Associated to the twisting map one has a three-point function $S\in 
C^\infty(M\times M\times M,\R)$ called 
the {\bf WKB-phase} of the elementary solvable symplectic symmetric space:
\begin{equation}
S(x_0,x_1,x_2):=\xi\left(\oint_{0,1,2}\sinh(a_0-a_1)l_2 \right);
\end{equation}
where $\oint_{0,1,2}$ stands for cyclic summation and where 
$x_i=(a_i,l_i)\quad(i=0,1,2)$. The phase $S$ turns 
out to be invariant under the (diagonal) action of the symmetries 
$\{s_x\}_{x\in M}$ on $M\times M\times M$.
This will be the essential constituent of the associative oscillatory
kernel defining a symmetry-invariant strict quantization on every
elementary solvable symplectic symmetric space.  We now recall this
construction as in \cite{Biep00a}.
\begin{dfn}
For a compactly supported function $u\in C^{\infty}_c(\p)$,
identifying $\l^\star$ with $\a$, we denote by $\tilde{u}\in C^\infty(\a\times\a)$ 
its partial Fourier transform:
\begin{equation}\label{FOURIER}
\tilde{u}(a,\alpha):=\int_\l
e^{i\Omega(\alpha,l)}u(a,l){\rm d}l.
\end{equation}
We also denote by $\phi_\hbar:\a\to\a$ the one-parameter family of 
twisting maps:
\begin{equation}
\phi_\hbar(a):=\frac{2}{\hbar}\phi(\frac{\hbar}{2}a).
\end{equation}
For $u,v\in C^{\infty}_c(\p)$, we set
\begin{equation}
<u\ |\ v>_\hbar :=\int_{\a\times \a}
\tilde{u}(a,\alpha)\overline{\tilde{v}(a,\alpha)}\,|\mbox{Jac}_{\phi^{-1}}(\alpha)|\,{\rm d}a\,{\rm d}\alpha.
\end{equation}
The pair $(C^\infty(\p),\,<\, ,\,>_\hbar)$ is a pre-Hilbert space, and
we denote by ${\cal H}_\hbar$ its Hilbert completion.
\end{dfn}
The Hilbert product $<\,,\,>_\hbar$ turns out to be symmetry-invariant on 
$C^\infty_c(M)$. The action of the transvection group then extends by 
continuity to an isometric action on ${\cal H}_\hbar$.

\begin{thm}\label{WKB}\cite{Biep00a} Let $(M,\omega,s)$ be a strictly
geodesically convex elementary solvable symplectic symmetric space. 
Realize it symplectically as $(\p=\a\times\l,\Omega)$, and define the
two-point function $A\in C^\infty(M\times M)$ by:
\begin{equation}\label{AMP}
A(x_1,x_2):=|\mbox{Jac}_{\phi}(a_1-a_2.)|
\end{equation}
This function is called the {\bf WKB-amplitude} and turns out to be 
symmetry-invariant.  In this notation, one has the following.
\begin{enumerate}
\item[(i)] For all $\hbar\in\R\backslash\{0\}$ and  $u,v\in C^\infty_c(M)$, the formula:

\begin{equation}\label{PROD}
u\star_\hbar v(x_0):=\int_{M\times M}u(x_1)
\,v(x_2)\,A(x_1,x_2)\,e^{\frac{i}{\hbar}S(x_0,x_1,x_2)}\,{\rm d}x_1\,{\rm d}x_2
\end{equation}

extends as an associative product on ${\cal H}_\hbar$ ($\,{\rm d}x$ 
denotes some normalization of the 
symplectic volume on $(M,\omega)$). Moreover, (for 
suitable $u,v$ and $x_{0}$) the stationary phase method yields a power 
series expansion of the form
\begin{equation}\label{STAR}
u\star_\hbar v(x_0)\sim 
uv(x_0)+\frac{\hbar}{2i}\{u,v\}(x_{0})+o(\hbar^{2});
\end{equation}
where $\{\,,\,\}$ denotes the symplectic Poisson bracket on $(M,\omega)$.
\item[(ii)] The pair $({\cal H}_\hbar,\star_\hbar)$ is a topological Hilbert 
algebra which the transvection group of $(M,\omega,s)$ acts on by 
automorphisms.
\end{enumerate}
\end{thm}
A classical procedure then produces a similar result in the 
$C^\star$-context, see \cite{Biep00a} for details.
\begin{rmk}
{\rm 
Wether a symmetric space is strictly geodesically convex is of course
entirely encoded in the spectral content of the splitting endomorphism
$\rho:\a\to\End(\b)$.  This is discussed in detail in \cite{Biep00a}. 
}
\end{rmk}

\section{Symplectic Lie groups associated to a class of symmetric
spaces}\label{SGS}
\begin{dfn}
A {\bf symplectic} Lie algebra is a pair $(\s,\omega)$ where $\s$
is a Lie algebra and $\omega\in\bigwedge^{2}(\s^{\star})$ is a non-degenerate 
Chevalley two-cocycle w.r.t. the trivial representation of $\s$.
\end{dfn}
In this section, we associate symplectic Lie algebras to a class
of symplectic symmetric spaces.
\subsection{Holonomy isotropic symplectic symmetric spaces}
\begin{dfn}
A symplectic triple $t=(\g,\sigma,\Omega)$ is called {\bf holonomy 
isotropic} ({\bf HI}) if $[\k,\p]$ is an isotropic subspace of 
$(\p,\Omega)$.
\end{dfn}
\begin{prop}\cite{Bi1}
A symplectic triple $t=(\g,\sigma,\Omega)$ is holonomy isotropic if and 
only if $[\g,\g]$ is Abelian.
\end{prop}
\begin{dfn}\label{SPLIT} Let $t=(\g,\sigma,\Omega)$ be HI and consider the
extension sequence
\begin{equation}
0\to[\g,\g]\to\g\to\a:=\g/[\g,\g]\to0.
\end{equation}
The HI triple $t$ is called {\bf split} if this extension is split.
\end{dfn}
\begin{lem}\label{HIS}
Let $t=(\g,\sigma,\Omega)$ be HI split. Set $\b=[\g,\g]$ and denote by 
$\rho:\a\to\End(\b)$ the splitting homomorphism. Then, realizing $\g$ as 
the semi-direct product $\g=\b\times_\rho\a$, one can assume that $\a$ 
is stable under $\sigma$.
\end{lem}
\Pf
For $a\in\a\subset\g$, write $a=a_\k+a_\p$ according to the decomposition 
w.r.t. $\sigma$. Then for all $a,a'\in\a$, one has
$0=[a,a']=[a_\p,a'_\p] + b\quad b\in[\k,\p]$ since $\k$ is Abelian.
This yields $[a_\p,a'_\p]=0$.

Therefore, for $\mbox{pr}_\p:\g\to\p$ the projection parallel to $\k$,
the $\p$-component $\mbox{pr}_\p(\a)$ is an Abelian
subalgebra of $\g$ supplementary to $\b$.  A dimension count then
yields the lemma.  \EPf

\begin{lem}
Assume that $t=(\g,\sigma,\Omega)$ is HI split, indecomposable and non-flat.
Set $0\to\b\to\g\to\a\to0$ as in Lemma~\ref{HIS}. Then $\a$ and $\l=[\k,\p]$
are in duality. In particular, there exists a $\k$-invariant symplectic structure 
on $\p$ for which $\a$ is Lagrangian.
\end{lem}

\Pf
Set $V:=\l^\perp\cap\a$ and choose a subspace $W$ of $\a$ in duality with 
$\l$. Then counting dimensions yields $\a=W\oplus V$. Moreover, in the 
decomposition $\p=\l\oplus W\oplus V$, the matrix of $\Omega$ is of the form
$$[\Omega]=
\left(
\begin{array}{ccc}
0 & I & 0\\
-I & 0 &B\\
0&-B'&A
\end{array}
\right).
$$
Since $\det[\Omega]\neq0$, one gets $\det\left(\begin{array}{cc} -I &B\\ 
0 &A\end{array} \right)\neq0$; hence $\det A\neq0$ and $V$ is symplectic.
Now, $\Omega([\k,V],\p)=\Omega(V,\l)=0$, hence $[\k,V]=0$. Also $[V,\l]
=[V,[\k,\p]]=0$ by Jacobi. Thus $V$ is central, and therefore trivial by 
indecomposability.
\EPf

\noindent We now assume that $(\g^1,\sigma^1)$ is the involutive Lie
algebra underlying a split HI symplectic triple which is
indecomposable and non-flat.  We fix $\Omega^1$ such that the HI
symplectic triple $t^1=(\g^1,\sigma^1,\Omega^1)$ with
$0\to\b^1\to\g^1\to\a^1\to0$ has $\a^{1}$ and $\l^1=[\k^1,\p^1]$
dual Lagrangian subspaces.  We then consider the associated exact triple that
we denote by $t=(\g,\sigma,\Omega)$ (if $t^1$ is already exact we set
$t=t^1$).  Observe that, since $\a^1$ is isotropic, the triple $t$ is
elementary solvable with
$$
0\to\b:=[\g,\g]\to\g\to\a:=\a^1\to0.
$$
\noindent We now follow the same procedure in \cite{Biep00a}.
The map $\rho:\a\to \End(\b)$ is injective (because $\Omega$ is
nondegenerate), so we may identify $\a$ with its image~:
$\a=\rho(\a)$.  Let $\Sigma:\End(\b)\to \End(\b)$ be the automorphism
induced by the conjugation with respect to the involution
$\sigma|_{\b}\in GL(\b)$, i.e. $\Sigma = Ad( \sigma|_{\b})$.  The
automorphism $\Sigma$ is involutive and preserves the canonical Levi
decomposition $\End(\b)={\cal Z}\oplus sl(\b)$, where ${\cal Z}$
denotes the center of $\End(\b)$.  Writing the element
$a=\rho(a)\in\a$ as $a=a_Z+a_0$ with respect to this decomposition,
one has~: $\Sigma(a)=a_Z+\Sigma(a_0)=-a=-a_Z-a_0$, because the
endomorphisms $a$ and $\sigma|_{\b}$ anticommute.  Hence
$\Sigma(a_0)=-2a_Z-a_0$ and therefore $a_Z=0$.  So, $\a$ actually lies
in the semisimple part $sl(\b)$.  For any $x\in sl(\b)$, we denote by
$x=x^S+x^N,\quad x^S,x^N\in sl(\b)$, its abstract Jordan-Chevalley
decomposition.  Observe that, for $sl(\b)=sl_+\oplus sl_-$, the
decomposition in $(\pm 1)$-$\Sigma$- eigenspaces, one has $\a\subset
sl_-$.  Also, $\a_N:=\{ a^N\}_{a\in\a}$ is an Abelian subalgebra in
$sl_-$ commuting with $\a$.  Set $\a_S:=\{a^S\}_{a\in\a}$.

\noindent Consider the complexification $\b^{c}:=\b\otimes\C$ and
$\C$-linearly extend the endomorphisms $\{\rho(a)\}_{a\in\a}$ and
$\sigma$.  Also consider the complex Lie algebra
$sl(\b^c)=sl(\b)\otimes\C$ and $\C$-linearly extend to $sl(\b^c)$ the
involution $\Sigma$.

\noindent Let
\begin{equation}
\b^c=:\bigoplus_{\alpha\in\Phi}\b_\alpha
\end{equation}
be the weight space decomposition w.r.t. the action of $\a_S$. Note that
for all $\alpha$, one has $\a_N.\b_\alpha\subset\b_\alpha$. Moreover,
for all $X_\alpha\in\b_\alpha$ and $a^S\in\a_S$, one has
\begin{eqnarray*}
\sigma(a^S.X_\alpha)=\alpha(a^S)\sigma(X_\alpha)
=\sigma a^S\sigma^{-1}\sigma 
X_\alpha=\Sigma(a^S).\sigma(X_\alpha)=-a^S.\sigma(X_\alpha).
\end{eqnarray*}
Therefore, $-\alpha\in\Phi$ and $\sigma\b_\alpha=\b_{-\alpha}$.
Note in particular that $\sigma\b_0=\b_0$.

\begin{lem}\label{BO}
If the triple $t^{1}$ is assumed indecomposable and non-flat, then
$$
\b_0=0.
$$
\end{lem}
\Pf Assume $0\in\Phi$.
For all $\alpha\in\Phi$, the subspace 
$$
V_\alpha:=\b_\alpha\oplus\b_{-\alpha}
$$
of $\b^c$ is stable under $\sigma$. In particular, the complexified 
involutive Lie algebra $(\g^c:=\g\otimes\C,\sigma)$ can be 
expressed as $\g^c=\a^c\times_{\rho}\b^c$ with
$$
\b^c=\bigoplus_{\alpha\in\Phi^+}\b_\alpha\oplus\b_0,
$$
where the `positive system' of weights $\Phi^+$ is chosen so that
$$
\Phi=\{0\}\cup\Phi^+\cup(-\Phi^+)\mbox{ (disjoint union.)}
$$
One therefore has the decomposition
$$
V_\alpha=\k_\alpha\oplus\l_\alpha
$$
into ($\pm$)-eigenspaces for $\sigma$. Moreover, since 
$\g^c=[\g^c,\g^c]$ and $[\a,\b_\alpha]\subset\b_\alpha$, one has:
\begin{eqnarray*}
\l_\alpha=[\a,\k_\alpha]\mbox{ and}\\
\k_\alpha=[\a,\l_\alpha],
\end{eqnarray*}
for all $\alpha\in\Phi^+\cup\{0\}$. This implies
$\l_0=[\a_N,\k_0]$ and $\k_0=[\a_N,\l_0]$. Hence
$\l_0=[\a_N,[\a_N,\l_0]]$ and an induction yields $\l_0=0$.
\EPf
\begin{cor}
A nilpotent HI split symplectic symmetric space is flat.
\end{cor}
\begin{prop}\label{SG}
Let  $t=(\g=\b\times_\rho\a,\sigma,\Omega=\delta\xi)$ be the exact triple
associated with a non-flat indecomposable split symplectic triple $t^1$.
Let $\Phi$ be the set of weights associated with the (complex) action 
of $\a_S$ on $\b^c$. Fix a positive system $\Phi^+$ and set
$$
\b^+:=\bigoplus_{\alpha\in\Phi^+}\b_\alpha.
$$
Then the pair $(\s^c:=\a^c\times_\rho\b^+,\Omega|_{\s^c})$ is a (complex)
symplectic Lie algebra.
\end{prop}
\Pf
By the proof of Lemma~\ref{BO}, the restricted projection
$\b^+\stackrel{p}{\to}\l^c:X\mapsto\frac{1}{2}(X-\sigma(X))$ is a linear isomorphism. 
Moreover, for all $X\in\b^+,a\in\a^c$, one has $\Omega(X,a)=\xi[p(c),a]$.
The proposition follows from the non-degeneracy of the pairing
$\a^c\times\l^c\to\C$.  \EPf

\begin{dfn}
Let $t$ be a HI split symplectic triple.
Decompose $t$ into a direct sum of indecomposables and a flat factor.
Proposition~\ref{SG} then canonically associates to $t$ a (complex)
symplectic Lie algebra $\s^c(t)$, the {\bf complex
symplectic Lie algebra associated with} $t$.
\end{dfn}
Combining the results of the present section with Section \ref{SQ}, one gets
\begin{thm} Let $\s$ be a symplectic Lie algebra associated with the
exact triple of a strictly geodesically convex HI split symplectic
symmetric space $M$.  Then $M$ is the manifold underlying the
connected simply connected Lie group whose Lie algebra is $\s$. 
Moreover, Theorem \ref{WKB} defines a left-invariant strict
deformation quantization of this symplectic Lie group.
\end{thm}

\begin{rmk}
{\rm 
\begin{enumerate}
\item Passing to a formal star product by the a stationary phase
method, one gets a left-invariant star product on every (group) direct
factor of this symplectic Lie group.  \item In the case of the
non-Abelian two-dimensional Lie algebra, strict universal deformation
formulae associated with non-oscillatory kernels were studied in
\cite{BiMae}.
\end{enumerate}
}
\end{rmk}

\section{Elementary solvable pre-symplectic Lie groups and symmetric spaces}
\subsection{A class of solvable Lie groups}
\begin{dfn}
A Lie group is called {\bf pre-symplectic} if it carries a 
left-invariant Poisson structure. 
\end{dfn}
One then observes
\begin{prop}
Let $(G,\pw)$ be a pre-symplectic Lie group with neutral element $e$. Then,
\begin{enumerate}
\item[(i)]the orthodual $\s$ of the radical of $\pw_{e}$,
$$
\s:=(\mbox{rad }\pw_e)^\perp,
$$
is a Lie subalgebra of the Lie algebra $\g$ of $G$;
\item[(ii)] the symplectic leaves of $\pw$ are the left classes of the 
analytic subgroup $S$ of $G$ whose Lie algebra is $\s$.
\end{enumerate}
\end{prop}
\noindent In particular, the Lie group $S$ is symplectic in the sense of 
Lichn\'erowicz. 
Symplectic Lie groups often tend to be solvable \cite{Li}.

\begin{dfn}
A symplectic Lie algebra $(\s,\omega)$ is called {\bf elementary 
solvable} if
\begin{enumerate}
\item[(i)] it is a split extension of Abelian Lie algebras $\a$ and $\d$:
\begin{equation}
0\longrightarrow\d\longrightarrow\s\longrightarrow\a\longrightarrow 
0;
\end{equation}
\item[(ii)] The cocycle $\omega$ is exact.
\end{enumerate}
\end{dfn}

\begin{prop}
Every elementary solvable symplectic Lie algebra is associated with a split HI 
symplectic symmetric space.
\end{prop}
\Pf
Denote by $\rho:\a\to\End(\d)$ the splitting homomorphism and by
$\overline{\rho}:\a\to\End(\d)$ the opposite 
representation: $\overline{\rho}(a)(X):=-\rho(a)(X),\quad X\in\d$.
Set 
$$
\b:=\d\oplus\d
$$
and let $\a$ act on $\b$ via
$\rho\oplus\overline{\rho}$. Define the involution $\sigma_\b$ of $\b$
by
$$
\sigma_\b(X,Y)=(Y,X),\qquad X,Y\in\d.
$$
Set 
$$
\g:=\b\times_{\rho\oplus\overline{\rho}}\a
$$
and define the involution $\sigma$ of $\g$ as
$$
\sigma:=\sigma_\b\oplus(-\mbox{id}_\a).
$$
One then observes that $(\g,\sigma)$ is an involutive Lie algebra. Note 
that $\k=\{(X,X)\}_{X\in\d}$ while $\p=\{(X,-X)\}_{X\in\d}$.

\noindent Let $\eta\in\d^\star$ be such that $\delta\eta=\omega$ and define
$\xi\in\k^\star$ by
$$
\xi(X,X):=\eta(X),\qquad X\in\d.
$$
Extending $\xi$ to $\g$ by $0$ on $\p$, one defines a symplectic coboundary
on $\g$:
$$
\Omega:=\delta\xi.
$$
The triple $(\g,\sigma,\Omega)$ then defines the desired elementary solvable
symplectic symmetric space.
\EPf
\begin{dfn}\label{ELE} Let $(\s,\omega)$ be an elementary solvable symplectic
Lie algebra.  Consider a split Abelian extension of $\s$:
\begin{equation}
0\longrightarrow\q\longrightarrow\g\longrightarrow\s\longrightarrow 
0.
\end{equation}
Then $\g$ has a canonical a Poisson structure $\pw_{e}$ whose associated
symplectic Lie algebra is $(\s,\omega)$ (for $q_i^*\in\q^\star$ and
$s_i^*\in\s^\star\quad(i=1,2)$, one has
$\pw_{e}((q_1^*,s^*_1),(q_2^*,s^*_2)):=\omega({}^\sharp
s_1^*,{}^\sharp s^*_2)$ where $\sharp:\s^\star\to\s$ denotes the
musical isomorphism).

\noindent Such a (pre-symplectic) split Abelian extension of an elementary
symplectic Lie algebra is called an {\bf elementary solvable pre-symplectic} Lie 
algebra.
\end{dfn}

\section{Universal Deformation Formulae}\label{SUDF} 
\noindent Let
$(\g=\q\times\s,\pw_e)$ be an elementary solvable pre-symplectic Lie
algebra with associated symplectic Lie algebra $\s$.  Consider the
associated connected simply connected Lie groups $G,Q$ and $S$, and
define the chart:
\begin{equation}\label{CHART}
Q\times S\rightarrow G:(q,s)\mapsto qs.
\end{equation}
{\bf We assume this is a global diffeomorphism}. 

\begin{thm}\label{INDUC}
Let $({\bf H}^S,\star^S)$ be an associative algebra of functions on $S$ such 
that ${\bf H}^S\subset \mbox{Fun}(S)$ is an invariant subspace w.r.t. the 
left regular representation of $S$ on $\mbox{Fun}(S)$ and $\star^S$ is a
$S$-left-invariant product on ${\bf H}^S$. Then the function space 
${\bf H}:=\{u\in\mbox{Fun}(G)\ |\ \forall p\in Q:u(q,.)\in {\bf H}^S\}$ is an 
invariant subspace of $\mbox{Fun}(G)$ w.r.t. the left regular representation 
of $G$ and the formula
\begin{equation}
u\star v(q,s):=(u(q,.)\star^Sv(q,.))(s)
\end{equation}
defines a $G$-left-invariant associative product on ${\bf H}$. 
\end{thm}
\Pf
The only thing to check is the left-invariance.
We begin by writing $\s$ as $\s=\b\times\a$ as in Section \ref{SGS}. 
Let $A$ and $B$ be the corresponding subgroups in $S$.  Within the
chart (\ref{CHART}), the group multiplication reads as follows. 
First, for $q,q'\in Q,s'\in S$, one has $q.(q',s')=(qq',s')$. 
Moreover for $s=ab\in S$, one has
$s.(q',s')=sq's'=abq's'=aq'a^{-1}abs'$ because $[\b,\q]=0$.  Hence
$s.(q',s')=(aq'a^{-1},ss')$.  This immediately implies the first
assertion .  Moreover, one has
\begin{eqnarray*}
(L_q^\star u\star L_q^\star v)(q',s')
&=&((L_q^\star u(q',.)\star^S(L_q^\star v(q',.))(s')
=(u(qq',.)\star^Sv(qq',.))(s')
= L_q^\star(u\star v)(q',s')\mbox{ and}\\
(L_s^\star u\star L_s^\star v)(q',s')
&=&(((L_s^\star u)(q',.)\star^S((L_s^\star v)(q',.))(s')
= L_s^\star (u(aq'a^{-1},.))\star^SL_s^\star (v(aq'a^{-1},.))(s')\\
&=&L_s^\star(u(aq'a^{-1},.)\star^Sv(aq'a^{-1},.))(s')
= u\star v(aq'a^{-1},ss')
= L_s^{\star}(u\star v)(q',s').
\end{eqnarray*}
\EPf
Of course, this also holds at the formal level i.e. in the case 
where ${\bf H}^{S}=C^\infty(S)[[\nu]]$ and $\star^S=\star^S_\nu$
is a left-invariant formal star product. Moreover, one has
\begin{cor}
In the formal case, the $G$-invariant star product on $G$ defined in 
Theorem \ref{INDUC} deforms the usual pointwise product of functions in the
direction of the left-invariant Poisson structure $\tilde{w}$, provided
$\star^S_\nu$ does so on $S$.
\end{cor}

\section{Crossed, Smash and Co-products}

Every algebra, coalgebra, bialgebra, Hopf algebra and vector space is 
taken over the field $k=\mathbb{R}$ or $\mathbb{C}$.  For classical
definitions and facts on these subjects, we refer to \cite{Swem69a}, 
\cite{Abee80a} or more fundamentally to \cite{MiMo65}.

\noindent To calculate with a coproduct $\Delta$, we use the Sweedler notation
\cite{Swem69a}: 
\begin{equation}
\Delta(b)=\sum_{(b)} b_{(1)}\otimes b_{(2)}.
\end{equation} 
With this notation, we can write coassociativity 
$\ds \Delta^{(2)}(b) := (\Delta\otimes \id )\circ \Delta (b) = (\id \otimes
\Delta )\circ \Delta (b)$ in the following way:
\begin{equation}
\Delta^{(2)}(b) := \sum_{(b)} b_{(1)}\otimes b_{(2)}\otimes b_{(3)}
= \sum_{(b)(b_1)} b_{(1)(1)}\otimes b_{(1)(2)}\otimes b_{(2)}
= \sum_{(b)(b_2)} b_{(1)}\otimes b_{(2)(1)}\otimes b_{(2)(2)}.
\end{equation}
By associativity one can define 
$\ds \Delta^{(n)}(b):= \sum_{(b)} b_{(1)}\otimes b_{(2)}\otimes
\cdots \otimes b_{(n+1)}\ ,\ n\in \mathbb{N}.
$
Cocommutativity means 
\begin{equation}
\Delta^{(n)}(b)=\sum_{(b)} b_{(\sigma (1))}\otimes b_{(\sigma(2))}\otimes
\cdots \otimes b_{(\sigma(n+1))}\ ,\ \forall\sigma\in \mathfrak{S}_n.
\end{equation}

We need two classical definitions \cite{Abee80a}:

\begin{dfn}\label{modalg}
Let $(B, ., \Delta_B)$ be a bialgebra and $(C, .)$ an algebra. 
$C$ is a (left-) {\bf $B$-module algebra} if $C$ is a $B$-module
and if the product on $C$ is a $B$-module map w.r.t. the canonical
$B$-module structure on $C\otimes C$ defined by the coproduct of $B$.
This means, for all $a,b \in B$ and $f,g\in C$:
\begin{enumerate}
\item $\ds (ab)\rightharpoonup f = a\rightharpoonup(b\rightharpoonup f)$ 
\item  $\ds a\rightharpoonup(f.g) 
= \sum_{(a)} (a_{(1)}\rightharpoonup f). 
(a_{(2)}\rightharpoonup g)
= \mu_C (\Delta_B(a)\rightharpoonup (f\otimes g))$
where $\mu_C$ is the product on $\ C$.
\end{enumerate} 
\end{dfn}

\begin{dfn}\label{modcoalg}
Let $(B, ., \Delta_B)$ be a bialgebra and $(D, \Delta_D)$ a coalgebra. 
$D$ is a (left-) {\bf $B$-module coalgebra} if $D$ is a $B$-module
and if the coproduct $\Delta_D$ on $D$ is a $B$-module map w.r.t. the canonical
$B$-module structure on $D\otimes D$ defined by the coproduct of $B$.
This means, for all $a,b \in B$ and $f\in D$:
\begin{enumerate}
\item $\ds (ab)\rightharpoonup f = a\rightharpoonup(b\rightharpoonup f) $ 
\item  $\ds \Delta_D(a\rightharpoonup f) 
= \Delta_B(a)\rightharpoonup \Delta_D(f)
= \sum_{(a)(f)} (a_{(1)}\rightharpoonup f_{(1)})\otimes (a_{(2)}
\rightharpoonup f_{(2)}).$
\end{enumerate} 
\end{dfn}

\noindent Combining these two definitions, we get

\begin{dfn} Let $(B, ., \Delta_B)$ and $(C, . , \Delta_C)$  be two bialgebras. 
$C$ is a (left-) {\bf $B$-module bialgebra} if $C$ is both a $B$-module algebra 
and a $B$-module coalgebra.
\end{dfn}

\noindent Finally we recall the following classical definition \cite{McL75a}:

\begin{dfn} 
Let $(B, .)$ be a algebra and $C$ a vector space. $C$ is a {\bf $B$-bimodule}
(or left-right $B$-module) if $C$ is a left $B$-module and a right $B$-module
with the following compatibility condition:
\begin{equation} 
(a\rightharpoonup f)\leftharpoonup b = a\rightharpoonup (f\leftharpoonup b)\ ,
\ \forall a,b\in B, \forall f\in C.
\end{equation}
\end{dfn}

\subsection{Crossed products}

The literature on Hopf algebras contains a large collection of what we
can call generically semi-direct products, or crossed products. Let us describe 
some of them. The simplest example of these crossed products is usually called 
the smash product (see \cite{Swem68a, Mol77a}):

\begin{dfn}\label{smash}
Let $B$ be a bialgebra and $C$ a $B$-module algebra. The smash product
$C\sharp B$ is the algebra constructed on the vector space $C\otimes B$ 
where the multiplication  is defined by
\begin{equation}
(f\otimes a)\stackrel{\rightharpoonup}{\star}(g\otimes b)
=\sum_{(a)}f\ (a_{(1)}\rightharpoonup g)\otimes a_{(2)} b
\end{equation}
for $f,g\in C$ and $a,b\in B$.
\end{dfn}

\begin{rmk}\hfill
{\rm
\begin{enumerate}
\item Verifying associativity is a direct calculation (cf. \ref{LRsmash}).
\item The smash product can be seen as the algebraic version
of what is called ``crossed product''
in the $C^*$-algebra literature \cite{DVDZ99a, Pedg79a}.
\item 
\begin{enumerate}
\item Let $H$ and $K$ be groups and let $\tau: K\rightarrow \auto (H) $ be
an action of $K$ on $H$. This induces a $k[K]$-module algebra structure
on $k[H]$. Then $ k[H]\sharp  k[K] \cong k[H \rtimes K]$,
$H \rtimes K$ denoting the semi-direct product of $H$ by $K$.
\item Similarly, for Lie algebras $\mathfrak{h}$ and $\mathfrak{k}$,
a Lie algebra homomorphism $\sigma:\mathfrak{k}\rightarrow \der(\mathfrak{h})$ 
induces a $\mathsf{U}\mathfrak{k}$-module algebra structure on 
$\mathsf{U}\mathfrak{h}$ ($\mathsf{U}\mathfrak{h}$, resp.
$\mathfrak{k}$, denoting the
universal envelopping algebra of $\mathfrak{h}$, resp. $\mathfrak{g}$). 
Then $\mathsf{U}\mathfrak{h}\sharp\mathsf{U}\mathfrak{k}\cong 
\mathsf{U}(\mathfrak{h}\rtimes \mathfrak{k})$.
\end{enumerate}

\item This product can be seen in the cohomological interpretation of Sweedler
\cite{Swem68a} as a representative of the trivial class of a theory of 
extensions. The formula of the smash product can be ``twisted" a little more
by some 2-cocycle from $B\otimes B$ to $C$ and is called a crossed product.
\item If $B$ and $C$ are bialgebras, $C$ a $B$-module algebra 
and $B$ a $C$-module algebra, with some compatibilities between the
two actions, one can write some kind of more ``symmetric" formula. 
S. Majid has called double crossproduct the resulting algebra \cite{Majs90a}.
\item If $C$ is a bialgebra and $B$ is cocommutative, the natural tensor 
coproduct on $C\otimes B$ yields a bialgebra structure on $C\sharp B$. If
everything is Hopf, $C\sharp B$ can be made Hopf as well \cite{Mol77a}.
\item By dualizing Definition \ref{smash}, one gets a coalgebra
called the cosmash product. Combining smash and cosmash in order to form
a bialgebra leads to the notion of bicrossproduct \cite{Majs90a}.
\end{enumerate}
}
\end{rmk}

\noindent Assuming $B$ is cocommutative we now introduce a generalization of the 
smash product.
\begin{dfn}\label{LRsmash}
Let $B$ be a cocommutative bialgebra and $C$ a $B$-bimodule algebra 
(i.e. a $B$-module algebra for both, left and right, $B$-module structures). 
The {\bf L-R-smash product}
$C\natural B$ is the algebra constructed on the vector space $C\otimes B$ 
where the multiplication is defined by
\begin{equation}
(f\otimes a)\star (g\otimes b)=\sum_{(a)(b)}
                (f\leftharpoonup b_{(1)}) (a_{(1)}\rightharpoonup g)
                                        \otimes a_{(2)} b_{(2)}
\end{equation}
for $f,g\in C$ and $a,b\in B$.
\end{dfn}

\begin{rmk} 
{\rm The L-R-smash product $C\natural B$ is an associative algebra. 
Indeed, this is an easy adaptation of the proof of the associativity
for the smash product: we first compute\\%%%%%%%%%%%
$\ds\left( (f\otimes a)\star (g\otimes b)\right)\star (h\otimes c)
=\left(\sum_{(a)(b)}(f\ars b_{(1)})(a_{(1)}\als g)\otimes a_{(2)}b_{(2)}\right)
	\star (h\star c)$\\
$\ds =\sum_{(a)(b)(c)(a_{(2)}b_{(2)})}
	\left(\left((f\ars b_{(1)})(a_{(1)}\als g)\right)\ars c_{(1)}\right)
	\left( (a_{(2)}b_{(2)})_{(1)}\als h\right)
	\otimes (a_{(2)}b_{(2)})_{(2)}c_{(2)}$\\
$\ds =\sum_{(a)(b)(c)(c_{(1)})(a_{(2)})(b_{(2)})}
	\left( (f\ars b_{(1)})\ars c_{(1)(1)}\right)
	\left( (a_{(1)}\als g)\ars c_{(1)(2)}\right)
	\left(a_{(2)(1)}b_{(2)(1)}\als h\right)
	\otimes a_{(2)(2)}b_{(2)(2)})c_{(2)}$\\
$\ds =\sum_{(a)(b)(c)}
	\left( f\ars b_{(1)}c_{(1)}\right)
	\left( (a_{(1)}\als g)\ars c_{(2)}\right)
	\left(a_{(2)}b_{(2)}\als h\right)
	\otimes a_{(3)}b_{(3)}c_{(3)}.$\\
Now we compute\\ %%%%%%%%%%%%%
$\ds (f\otimes a)\star\left( (g\otimes b)\star (h\otimes c)\right)
=(f\star a)\star\left(\sum_{(b)(c)}(g\ars c_{(1)})(b_{(1)}\als h)
	\otimes b_{(2)}c_{(2)}\right)$\\
$\ds =\sum_{(a)(b)(c)(b_{(2)}c_{(2)})}
	\left(f\ars (b_{(2)}c_{(2)})_{(1)}\right)
	\left(a_{(1)}\als \left((g\ars c_{(1)})(b_{(1)}\als h)\right)\right)
	\otimes a_{(2)}(b_{(2)}c_{(2)})_{(2)}$\\
$\ds =\sum_{(a)(a_{(1)})(b)(c)(b_{(2)})(c_{(2)})}
	\left(f\ars b_{(2)(1)}c_{(2)(1)}\right)
	\left(a_{(1)(1)}\als (g\ars c_{(1)})\right)
	\left(a_{(1)(2)}\als (b_{(1)}\als h)\right)
	\otimes a_{(2)}b_{(2)(2)}c_{(2)(2)}$\\	
$\ds =\sum_{(a)(b)(c)}
	\left( f\ars b_{(2)}c_{(2)}\right)
	\left( (a_{(1)}\als g)\ars c_{(1)}\right)
	\left(a_{(2)}b_{(1)}\als h\right)
	\otimes a_{(3)}b_{(3)}c_{(3)}.$\\
Since $B$ is cocommutative we have the result.
}
\end{rmk} 
In the same spirit, one has
\begin{lem}\label{cop}
If $C$ is a $B$-bimodule bialgebra, the natural tensor 
product coalgebra structure on $C\otimes B$ defines a bialgebra 
structure to $C\natural B$. 

If $C$ and $B$ are Hopf algebras, $C\natural B$ is a Hopf %%%%%%%
algebra as well, defining the antipode by
\begin{eqnarray}
J_\star(f\otimes a) & = & \sum_{(a)}J_B(a_{(1)})\als J_C(f)\ars J_B(a_{(2)})
	\otimes J_B(a_{(3)})\\ 
	& = & \sum_{(a)}(1_C\otimes J_B(a_{(1)}))\star (J_C(f)\otimes 1_B)\star 
	(1_C\otimes J_B(a_{(2)})).\nonumber
\end{eqnarray}
\end{lem} 
\Pf The same type of calculation as above. \EPf

\noindent Under certain conditions, the L-R-smash product can be decomposed in 
``smash product like" terms:

\begin{lem}
\begin{enumerate} 
\item[(i)] If $B$ is commutative, one has
 \begin{eqnarray*}
(f\otimes a)\star (g\otimes b) 
&=&\left((f\otimes 1) \stackrel{\leftharpoonup}{\star} (1\otimes b)\right)
.\left((1\otimes a)\stackrel{\rightharpoonup}{\star}(g\otimes 1)\right)\\
&=&  \left( \sum_{(b)}(f \leftharpoonup b_{(1)})\otimes b_{(2)}\right)
.  \left(\sum_{(a)}(a_{(1)}\rightharpoonup g)\otimes a_{(2)}\right).
  \end{eqnarray*}

\item[(ii)] If $C$ is commutative, one has
  \begin{eqnarray*}
    (f\otimes a)\star (g\otimes b)
&=& \left((1\otimes a)\stackrel{\rightharpoonup}{\star}(g\otimes 1)\right)
. \left((f\otimes 1) \stackrel{\leftharpoonup}{\star} (1\otimes b)\right)\\
&=&  \left(\sum_{(a)}(a_{(1)}\rightharpoonup g)\otimes a_{(2)}\right)
.  \left( \sum_{(b)}(f \leftharpoonup b_{(1)})\otimes b_{(2)}\right).
  \end{eqnarray*}
  
\end{enumerate} 
\end{lem}
\noindent The proof is straightforward.

\noindent Now by a careful computation, one proves
\begin{prop}\label{*S}
Let $B$ be a cocommutative bialgebra, $C$ a $B$-bimodule algebra
and $(C\natural B, \star )$ their L-R-smash product.\\ 
Let $S$ be a linear automorphism of $C$ (as a vector space). We define:
\begin{enumerate}
\item[(i)]  the product $\bullet^S$ on $C$ by 
\begin{equation}
f\bullet^S g=S^{-1}\left(S(f).S(g)\right);
\end{equation}
\item[(ii)]  the left and right $B$-module structures, $\als^S$ and 
$\ars^S$, by 
\begin{equation}
a\als^S f := S^{-1}\left(a\rightharpoonup S(f)\right)\ \mbox{ and }\
f\ars^S a := S^{-1}\left(S(f)\leftharpoonup a\right);
\end{equation}
\item[(iii)]  the product, $\star^S$, on $C\otimes B$ by 
\begin{equation}
(f\otimes a)\star^S (g\otimes b)
=T^{-1}\left( T(f\otimes a) \star T(g\otimes b)\right)
\end{equation}
where $T:=S\otimes \id $.
\end{enumerate} 
Then $(C, \bullet^S)$ is a $B$-bimodule algebra for $\als^S$ and $\ars^S$ 
and $\star^S$ is the L-R-smash product defined by these structures.

\noindent Moreover, 
if $(C, . , \Delta_C , J_C , \rightharpoonup , \leftharpoonup )$ 
is a Hopf algebra  and a $B$-module bialgebra, then 
$$C_S:=(C, \bullet^S , \Delta^S_C:=(S^{-1}\otimes S^{-1})\circ\Delta_C\circ S ,
J^S_C:=S^{-1}\circ J_C\circ S , \als^S , \ars^S )$$ 
is also a Hopf algebra  and a $B$-module bialgebra. 
Therefore, by Lemma \ref{cop}, 
$$(C_S\natural B , \star^S , \Delta^S = (23)\circ (\Delta^S_C\otimes 
\Delta_B) , J_\star^S ),$$
is a Hopf algebra for
$\Delta^S$ the natural tensor product coalgebra structure
on $C_S\natural B$ (with $(23):C\otimes C\otimes B\otimes B\rightarrow 
C\otimes B\otimes C\otimes B,\
c_1\otimes c_2\otimes b_1\otimes b_2\mapsto
c_1\otimes b_1\otimes c_2\otimes b_2$)
and $J_\star^S$ the antipode given on $C_S\natural B$ by Lemma \ref{cop}.
Also, one has $$\Delta^S=(T^{-1}\otimes T^{-1})\circ (23)
\circ (\Delta_C\otimes\Delta_B)\circ T\ \mbox{ and }\
J_\star^S=T^{-1}\circ J_\star \circ T$$ with $T=S\otimes \id $.
\end{prop}

\section{Examples in deformation quantization}\label{EX}

\subsection{A construction on $T^\star(G)$}
Let  $G$ be a Lie group with Lie algebra $\g$ and $T^\star(G)$ its cotangent 
bundle. We denote by $\mathsf{U}\mathfrak{g}$,
$\mathsf{T}\mathfrak{g}$ and $\mathsf{S}\mathfrak{g}$ respectively the 
enveloping, tensor and symmetric algebras of $\mathfrak{g}$. 
Let $\mathsf{Pol}(\mathfrak{g}^*)$ be
the algebra of polynomial functions on $\mathfrak{g}^*$. 
We have the usual identifications:
$$\mathcal{C}^\infty (T^*G)\simeq \mathcal{C}^\infty (G\times \mathfrak{g}^*)
\simeq \mathcal{C}^\infty(G)\hat{\otimes} 
\mathcal{C}^\infty(\mathfrak{g}^*) \supset 
\mathcal{C}^\infty(G)\otimes\mathsf{Pol}(\mathfrak{g}^*)
\simeq \mathcal{C}^\infty(G)\otimes\mathsf{S}\mathfrak{g}.$$

\noindent First we deform $\mathsf{S}\mathfrak{g}$ via 
the ``parametrized version'',$\mathsf{U}_t\mathfrak{g}$, 
of $\mathsf{U}\mathfrak{g}$ defined by
$$\mathsf{U}_t\mathfrak{g}= \frac{\mathsf{T}\mathfrak{g}[[t]]}
{<XY-YX-t[X,Y]; X,Y\in \mathfrak{g}>}.$$
$\mathsf{U}_t\mathfrak{g}$ is naturally a Hopf algebra with 
$\Delta(X)=1\otimes X + X\otimes 1$, $\epsilon(X)=0$ and 
$S(X)=-X$ for $X\in \mathfrak{g}$.
For $X\in\g$, we denote by $\widetilde{X}$ (resp. $\overline{X}$) the
left- (resp. right-) invariant vector field on $G$ such that 
$\widetilde{X}_e=\overline{X}_e=X$. We consider
the following $k[[t]]$-bilinear actions of $B=\mathsf{U}_t\mathfrak{g}$
on $C=\mathcal{C}^\infty(G)[[t]]$, for $f\in C$ and $\lambda\in [0,1]$:
\begin{enumerate}
\item[(i)] $(X \rightharpoonup f)(x) = t (\lambda -1)\ 
(\widetilde{X}.f)(x)$,
\item[(ii)] $(f \leftharpoonup X)(x) = t \lambda \ (\overline{X}.f)(x)$.
\end{enumerate}
One then has
\begin{lem}
$C$ is a $B$-bimodule algebra w.r.t. the above left and right actions 
(i) and (ii).
\end{lem}

\begin{dfn}
We denote by $\star_\lambda$ the star product on
$\left(\mathcal{C}^\infty(G)\otimes\mathsf{Pol}(\mathfrak{g}^*)\right)[[t]]$ 
given by the L-R-smash product on 
$\mathcal{C}^\infty(G)[[t]]\otimes\mathsf{U}_t\mathfrak{g}$
constructed from the bimodule structure of the preceding lemma.
\end{dfn}

\begin{prop}
For $G=\mathbb{R}^n$, $\star_{\frac{1}{2}}$ is the Moyal star product 
(Weyl ordered), $\star_{0}$ is the standard ordered star product and 
$\star_{1}$ the anti-standard ordered one. 
In general $\star_\lambda$ yields the 
$\lambda$-ordered quantization, within the 
notation of M. Pflaum \cite{Pflm99a}. 
\end{prop}
\Pf
\noindent Let $\{ q_1, q_2,\ldots , q_n\}$ be coordinates on $\mathbb{R}^n$ and 
$\{ p_1, p_2,\ldots , p_n\}$ dual coordinates.
For
$\mathbf{l}=(l_1, l_2,\ldots , l_n)$ and $\mathbf{r}=(r_1, r_2,\ldots , r_n)$ 
in $\mathbb{N}^n$, set
$|\mathbf{l}|=l_1 + l_2 + \ldots + l_n$ and
$\mathbf{p}^{\mathbf{r}}=p_1^{r_1}p_2^{r_2}\ldots p_n^{r_n}$. 
Define
$$\frac{\partial^{|\mathbf{l}|+|\mathbf{m}|}u}
{\partial\mathbf{q}^\mathbf{m}\partial\mathbf{p}^\mathbf{l}}
= \frac{\partial^{|\mathbf{l}|+|\mathbf{m}|}u}
{\partial q_1^{m_1}\partial q_2^{m_2}\cdots\partial q_n^{m_n}
\partial p_1^{l_1}\partial p_2^{l_2}\cdots\partial p_n^{l_n}}, 
\quad u\in \mathcal{C}^\infty(\mathbb{R}^n\times {\mathbb{R}^n}^*).$$
For $\mathbf{r}!:=r_1! r_2!\ldots r_n!$ and 
$\ds \cbi{\mathbf{l}}{\mathbf{r}} 
:= \frac{\mathbf{r}!}{\mathbf{l}! \mathbf{r-l}!}
= \cbi{l_1}{r_1} \cbi{l_2}{r_2}\ldots  \cbi{l_n}{r_n}$,
we obtain \\ $\ds \Delta (\mathbf{p}^\mathbf{r})=\sum_{\mathbf{l}
=\mathbf{0}}^\mathbf{r} \cbi{\mathbf{l}}{\mathbf{r}} 
\mathbf{p}^\mathbf{l} \otimes \mathbf{p}^{\mathbf{r}-\mathbf{l}}$, 
for the coproduct in $\mathsf{Pol}({\mathbb{R}^n}^*)
=\mathsf{S}\mathbb{R}^n$.

\noindent Now
$\ds (f\otimes \mathbf{p}^\mathbf{r})\star_\lambda 
		(g\otimes \mathbf{p}^\mathbf{s})
=\sum_{(\mathbf{p}^\mathbf{r})}\sum_{(\mathbf{p}^\mathbf{s})}
\left(f\ars \mathbf{p}^\mathbf{s}_{(1)}\right)
\left(\mathbf{p}^\mathbf{r}_{(1)}\als g \right)
 \otimes \mathbf{p}^\mathbf{r}_{(2)}\mathbf{p}^\mathbf{s}_{(2)}$

\begin{eqnarray*}
& = & \sum_{(\mathbf{p}^\mathbf{r})} \sum_{\mathbf{m}=\mathbf{0}}^{\mathbf{s}}
  \cbi{\mathbf{m}}{\mathbf{s}} \left(f\ars \mathbf{p}^\mathbf{m}\right) 
\left(\mathbf{p}^\mathbf{r}_{(1)}\als g \right)
 \otimes \mathbf{p}^\mathbf{r}_{(2)}\mathbf{p}^{\mathbf{s}-\mathbf{m}}\\
& = & \sum_{\mathbf{l}=\mathbf{0}}^{\mathbf{r}}\sum_{\mathbf{m}
	=\mathbf{0}}^{\mathbf{s}}
 \cbi{\mathbf{l}}{\mathbf{r}} \cbi{\mathbf{m}}{\mathbf{s}}
 \left(f\ars \mathbf{p}^\mathbf{m}\right) 
	\left(\mathbf{p}^\mathbf{l}\als g \right)
 \otimes \mathbf{p}^{\mathbf{r}-\mathbf{l}} \mathbf{p}^{\mathbf{s}-\mathbf{m}}\\
&  = & \sum_{\mathbf{l}=\mathbf{0}}^{\mathbf{r}} \sum_{\mathbf{m}
	=\mathbf{0}}^{\mathbf{s}} t^{|\mathbf{l}|+|\mathbf{m}|} 
\cbi{\mathbf{l}}{\mathbf{r}} \cbi{\mathbf{m}}{\mathbf{s}}
\lambda^{|\mathbf{m}|} (\lambda-1)^{|\mathbf{l}|}
\frac{\partial^{|\mathbf{m}|}f}{\partial\mathbf{q}^\mathbf{m}}
 \frac{\partial^{|\mathbf{l}|}g}{\partial\mathbf{q}^\mathbf{l}}
 \otimes \mathbf{p}^{\mathbf{r}+\mathbf{s}-(\mathbf{l}+\mathbf{m})} 
\end{eqnarray*}

\noindent The formula for the $\lambda$-ordered star product 
$\star_\lambda$ \cite{Pflm99a} is 
$$ u *_\lambda v = \sum_{\mathbf{l}+\mathbf{m}\geqslant \mathbf{0}}
\frac{t^{|\mathbf{l}|+|\mathbf{m}|}}{\mathbf{l}!\mathbf{m}!} 
\lambda^{|\mathbf{m}|} (\lambda-1)^{|\mathbf{l}|}
\frac{\partial^{|\mathbf{l}|+|\mathbf{m}|}u}
		{\partial\mathbf{q}^\mathbf{m}\partial\mathbf{p}^\mathbf{l}}
\frac{\partial^{|\mathbf{l}|+|\mathbf{m}|}v}
		{\partial\mathbf{q}^\mathbf{l}\partial\mathbf{p}^\mathbf{m}}.$$
Hence $\ds (f(\mathbf{q}) \mathbf{p}^\mathbf{r})\star_\lambda 
	(g(\mathbf{q}) \mathbf{p}^\mathbf{s})
=\sum_{\mathbf{l}+\mathbf{m}\geqslant \mathbf{0}}
\frac{t^{|\mathbf{l}|+|\mathbf{m}|}}{\mathbf{l}!\mathbf{m}!} 
\lambda^{|\mathbf{m}|} (\lambda-1)^{|\mathbf{l}|}
\frac{\partial^{|\mathbf{l}|+|\mathbf{m}|}f(\mathbf{q}) \mathbf{p}^\mathbf{r}}
{\partial\mathbf{q}^\mathbf{m}\partial\mathbf{p}^\mathbf{l}}
\frac{\partial^{|\mathbf{l}|+|\mathbf{m}|}g(\mathbf{q}) \mathbf{p}^\mathbf{s}}
{\partial\mathbf{q}^\mathbf{l}\partial\mathbf{p}^\mathbf{m}}$

\begin{eqnarray*}
& = & \sum_{\mathbf{l}+\mathbf{m}\geqslant \mathbf{0}}
\frac{t^{|\mathbf{l}|+|\mathbf{m}|}}{\mathbf{l}!\mathbf{m}!} 
\lambda^{|\mathbf{m}|} (\lambda-1)^{|\mathbf{l}|}
\frac{\partial^{|\mathbf{l}|}\mathbf{p}^\mathbf{r}}
		{\partial\mathbf{p}^\mathbf{l}}
\frac{\partial^{|\mathbf{m}|}\mathbf{p}^\mathbf{s}}
		{\partial\mathbf{p}^\mathbf{m}}
\frac{\partial^{|\mathbf{m}|}f}{\partial\mathbf{q}^\mathbf{m}}
\frac{\partial^{|\mathbf{l}|}g}{\partial\mathbf{q}^\mathbf{l}}\\
&  = & \sum_{\mathbf{l}=\mathbf{0}}^{\mathbf{r}} \sum_{\mathbf{m}
	=\mathbf{0}}^{\mathbf{s}} 
\frac{t^{|\mathbf{l}|+|\mathbf{m}|}}{\mathbf{l}!\mathbf{m}!} 
\frac{\mathbf{r}!}{(\mathbf{r}-\mathbf{l})!}\mathbf{p}^{\mathbf{r}-\mathbf{l}}
\frac{\mathbf{s}!}{(\mathbf{s}-\mathbf{m})!}\mathbf{p}^{\mathbf{s}-\mathbf{m}}
\lambda^{|\mathbf{m}|} (\lambda-1)^{|\mathbf{l}|}
\frac{\partial^{|\mathbf{m}|}f}{\partial\mathbf{q}^\mathbf{m}}
\frac{\partial^{|\mathbf{l}|}g}{\partial\mathbf{q}^\mathbf{l}}\\
&  = & \sum_{\mathbf{l}=\mathbf{0}}^{\mathbf{r}} \sum_{\mathbf{m}
=\mathbf{0}}^{\mathbf{s}} t^{|\mathbf{l}|+|\mathbf{m}|}
\cbi{\mathbf{l}}{\mathbf{r}} \cbi{\mathbf{m}}{\mathbf{s}} 
\lambda^{|\mathbf{m}|} (\lambda-1)^{|\mathbf{l}|}
  \frac{\partial^{|\mathbf{m}|}f}{\partial\mathbf{q}^\mathbf{m}}
 \frac{\partial^{|\mathbf{l}|}g}{\partial\mathbf{q}^\mathbf{l}}
 \mathbf{p}^{\mathbf{r}+\mathbf{s}-(\mathbf{l}+\mathbf{m})}.
\end{eqnarray*}
\EPf
\begin{rmk}
{\rm 
In the general case, it would be interesting to 
compare our $\lambda$-ordered L-R smash product with classical
constructions of star products on $T^\star(G)$ with Gutt's product
as one example \cite{Guts83a}.
}
\end{rmk}

\subsection{Hopf structures}\label{CP}
We have discussed (see Lemma \ref{cop})
the possibility of having a Hopf structure on $C\natural B$.  Let
us consider the particular case of $\mathcal{C}^\infty
(\mathbb{R}^n)[[t]] \natural\mathsf{U}_t\mathbb{R}^n =
\mathcal{C}^\infty (\mathbb{R}^n)[[t]] \natural\mathsf{S}\mathbb{R}^n$
($\mathbb{R}^n$ is commutative).  $\mathsf{S}\mathbb{R}^n$ is endowed
with its natural Hopf structure but we also need a Hopf
structure on $\mathcal{C}^\infty (\mathbb{R}^n)[[t]] =
\mathcal{C}^\infty (\mathbb{R}^n)\otimes\mathbb{R}[[t]].$ We will not
use the usual one. Our alternative structure is defined as follows.
\begin{dfn} 
We endow $\mathbb{R}[[t]]$ with the usual product, the co-product
$\Delta(P)(t_1,t_2):=P(t_1+t_2)$, the co-unit $\epsilon(P)=P(0)$
and the antipode $J(t)=-t$.
We consider the Hopf algebra
$(\mathcal{C}^\infty (\mathbb{R}^n), . , 
{\bf 1}, \Delta_C , \epsilon_C , J_C )$, with %where $.$ is the 
pointwise multiplication, the unit ${\bf 1}$ 
(the constant function of value $1$), the coproduct
$\Delta_C (f)(x,y)=f(x+y)$, the co-unit $\epsilon (f)=f(0)$
and the antipode $J_C(f)(x)=f(-x)$. 
The tensor product of these two Hopf algebras
then yields a  Hopf algebra denoted by
$$(\mathcal{C}^\infty (\mathbb{R}^n)[[t]], . , {\bf 1} ,
\Delta_t , \epsilon_t , J_t).$$ 
Note that $\Delta_t$ and $J_t$ are not linear in $t$.
We then define, on the L-R smash $\mathcal{C}^\infty (\mathbb{R}^n)[[t]]
\natural\mathsf{S}\mathbb{R}^n$,
$$\Delta_\star := (23)\circ (\Delta_t\otimes\Delta_B),\quad
\epsilon_\star := \epsilon_t\otimes\epsilon_B
\quad \mbox{ and }\quad J_\star\ \mbox{ as in Lemma \ref{cop} }.$$ 
\end{dfn}

\begin{prop}
($\mathcal{C}^\infty (\mathbb{R}^n)[[t]]\natural\mathsf{S}\mathbb{R}^n,
\star_\lambda , {\bf 1}\otimes 1 , \Delta_\star , \epsilon_\star , J_\star )$
is a Hopf algebra.
\end{prop}
\Pf
According to Lemma \ref{cop}, the only thing to show is
that the actions $\als$ and $\ars$ are coalgebra maps.

\noindent Let $\{ X_i; i=1, ..., n\}$ be a basis of $\mathbb{R}^n$.
Since $\mathbb{R}^n$ is commutative, left and right invariant
vector fields coincide so it is enough to show that $\als$
is a coalgebra map, that is, for $a\in\mathsf{S}\mathbb{R}^n$
and $\tilde{f}\in\mathcal{C}^\infty (\mathbb{R}^n)[[t]]$,
$$\Delta_t(a\als \tilde{f}) = \Delta_B(a)\als \Delta_t(\tilde{f}).$$
By additivity it suffices to check this for 
$\tilde{f}=t^nf,\ n\in \mathbb{N}, f\in\mathcal{C}^\infty (\mathbb{R}^n)$, 
and since
\begin{eqnarray*}
\Delta_t(X_ia\als \tilde{f})&=&\Delta_t(X_i\als (a\als \tilde{f}))
	=\Delta_B(X_i)\als\Delta_t(a\als \tilde{f})
	=\Delta_B(X_i)\als(\Delta_B(a)\als\Delta_t(\tilde{f}))\\
	&=&(\Delta(X_i)(\Delta_B(a))\als\Delta_t(\tilde{f})
	=\Delta(X_ia)\als\Delta_t(\tilde{f}),
\end{eqnarray*}
it suffices to check this on the $X_i$'s. We have
\begin{eqnarray*}
&&\Delta_t(X_i\als t^nf)(X,Y)
	=\Delta_t(t^{n+1}\frac{\partial f}{\partial x_i})(X,Y)
	=\Delta_t(t^{n+1})\Delta_t(\frac{\partial f}{\partial x_i})(X,Y)
	=(t_1+t_2)^{n+1}\frac{\partial f}{\partial x_i}(X+Y)\\
&&\mbox{ and }(\Delta_b(X_i)\als\Delta_t(t^nf))(X,Y)
	=(X_i\otimes 1 + 1\otimes X_i)\als \left( (t_1+t_2)^nf(X+Y)\right)\\
	&=&t_1(t_1+t_2)^n(\frac{\partial}{\partial x_i}\otimes 1).f(X+Y)
	+t_2(t_1+t_2)^n(1\otimes\frac{\partial}{\partial x_i}).f(X+Y)\\
	&=&t_1(t_1+t_2)^n\frac{\partial f}{\partial x_i}(X+Y)
	+t_2(t_1+t_2)^n\frac{\partial f}{\partial x_i}(X+Y)
	=(t_1+t_2)^{n+1}\frac{\partial f}{\partial x_i})(X+Y)
\end{eqnarray*}
\EPf

\begin{rmk}
{\rm 
The case $\lambda=\frac{1}{2}$ yields the usual Hopf structure on 
the enveloping algebra of the Heisenberg Lie algebra.
}
\end{rmk}

\subsection{Symplectic symmetric spaces}
On an elementary 
solvable symplectic symmetric space we have seen that there exists a global 
Darboux chart such that
$(M,\omega )\simeq (\mathfrak{p}=\mathfrak{l}\oplus \mathfrak{a},\Omega)$.
Thus one has 
$$
\mathcal{C}^\infty (M)\simeq\mathcal{C}^\infty (\mathfrak{p})
\simeq \mathcal{C}^\infty (\mathfrak{l}) \hat\otimes
\mathcal{C}^\infty (\mathfrak{a}) \buildrel\simeq_{a\simeq\mathfrak{l}^{\ast}}^{}
\mathcal{C}^\infty
(\mathfrak{l}) \hat\otimes \mathcal{C}^\infty (\mathfrak{l}^\ast)
\supset \mathcal{C}^\infty
(\mathfrak{l})\otimes Pol(\mathfrak{l}^\ast)
\buildrel\simeq_{\mathfrak{l}\hbox{ {\tiny{Abelian}}}}^{}
\mathcal{C}^\infty (\mathfrak{l})\otimes \mathsf{U}\mathfrak{l}.
$$
WKB-quantization on such a space turns out to be a L-R-smash 
product. Namely, one has
\begin{prop}
The formal version of the invariant WKB-quantization of an elementary 
solvable symplectic symmetric spaces defined in Theorem \ref{WKB} is 
a L-R-smash product of the form
$\star^S$ (cf. Proposition \ref{*S}).
\end{prop}
\Pf
Let ${\cal S}(\l)$ denote the Schwartz space of the vector space $\l$.
It is shown in \cite{Biep00a} that the map
$$
\begin{array}{ccc}
{\cal S}(\l) &\stackrel{S}{\longrightarrow}&{\cal S}(\l) \\
u&\mapsto&S(u):=F^{-1}\phi_\hbar^\star F(u)
\end{array}
$$
is a linear injection for all $\hbar\in\R$ ($F$ denotes the partial 
Fourier transform (\ref{FOURIER})).
An asymptotic expansion in a power series in $\hbar$ then yields
a formal equivalence, again denoted by $S$:
$$
S:=\id+o(\hbar):C^\infty(\l)[[\hbar]]\to C^\infty(\l)[[\hbar]].
$$
Carrying the Moyal star product on $(\p=\l\times\a,\Omega)$ by
$\ds \ T:=S\otimes \id \ $ yields a star product on $M\simeq\p$ 
which coincides with
the asymptotic expansion of the invariant WKB-product (\ref{PROD}).
\EPf
\noindent Subsection \ref{CP} then yields
\begin{cor}
The UDF's for elementary pre-symplectic Lie groups constructed in 
Section \ref{SUDF} admit compatible co-products and antipodes.
\end{cor}

\end{document}